\documentclass[12pt]{article}
\usepackage{amsmath,latexsym}

\newtheorem{theorem}{Theorem}
\newtheorem{lemma}[theorem]{Lemma}  
\newtheorem{corollary}[theorem]{Corollary}

\def\f2{\mathbf{F}_2}
\def\fq{\mathbf{F}_q}
\newcommand{\F}{\mathbf{F}}
\newcommand{\glnq}{\mbox{\rm GL}_{n}(q)}
\newcommand{\gl}{\mbox{\rm GL}}
\newcommand{\mnq}{\mbox{\rm M}_{n}(q)}
\newcommand{\cpl}{c_{\phi}(\lambda)}
\def\g{\gamma}

\newcommand{\qed}{\quad \hspace*{\fill} $\Box$ \par}
\newcommand{\proof}{\noindent \textbf{Proof} \,\,}
\newcommand{\ra}{\rightarrow}
\newcommand{\qbinom}[2]{\genfrac{(}{)}{0pt}{}{#1}{#2}_{\! q}}
\title{\textbf{Integer Sequences and Matrices Over Finite Fields}
}
\author{
    Kent E. Morrison \\ 
    Department of Mathematics \\
    California Polytechnic State University \\
    San Luis Obispo, CA 93407 \\
    \\
     \texttt{kmorriso@calpoly.edu} \\
}
       
\date{February 1, 2006}
\begin{document}
\maketitle  

\begin{abstract}
\noindent In this expository article we collect the integer sequences that count several different types of matrices over finite fields and provide references to the Online Encyclopedia of Integer Sequences (OEIS). Section 1 contains the sequences, their generating functions, and examples. Section 2 contains the proofs of the formulas for the coefficients and the generating functions of those sequences if the proofs are not easily available in the literature. The cycle index for matrices is an essential ingredient in most of the derivations.

\vspace{2mm}

\vspace{2mm}
 \end{abstract}

\begin{center} \textbf{Notation}  \end{center}
\begin{tabbing}
99999999\=99999999\=   \kill
\> $q$ \>  a prime power \\
\> $\fq$ \> field with $q$ elements\\
\> $\glnq$  \> group of $n \times n$ invertible matrices over $\fq$\\
\> $\mnq$ \> algebra of $n \times n$ matrices over $\fq$ \\
\> $\g_n(q)$  \> order of the group $\glnq$ \\
\> $\g_n$   \> order of $\glnq$ with $q$ understood \\
\> $\nu_d$  \> number of irreducible monic polynomials of degree $d$ over $\fq$
\end{tabbing}
\pagebreak
\begin{center}
\begin{tabular}{lcr}
\textbf{Sequences} & \textbf{Section} &\textbf{OEIS} \\
  \\
All matrices & 1.1 & A002416  \\
Invertible matrices & 1.2 & A002884 \\
Subspaces, $q$-binomial coefficients &  1.3 &A022166--188 \\
&& A006116--122 \\ 
&& A015195--217 \\
Splittings (direct sum decompositions) & 1.4\\ 
$q$-Stirling numbers, $q$-Bell numbers & 1.4 \\
Flags of subspaces & 1.5 & A005329, A069777\\
Linear binary codes & 1.6 & A022166, A076831\\
Matrices by rank & 1.7 \\
Linear derangements & 1.8, 2.2 & A002820 \\
Projective derangements & 1.8, 2.2\\
Diagonalizable matrices & 1.9, 2.3 \\
Projections & 1.10, 2.4 & A053846\\
Solutions of $A^k=I$ & 1.11, 2.5 & A053718, 722, 725\\
&& A053770--777\\
&& A053846--849 \\
&& A053851--857 \\
&& A053859--863 \\
Nilpotent matrices & 1.12 & A053763\\
Cyclic(regular) matrices & 1.13, 2.6 \\
Semi-simple matrices & 1.14, 2.7 \\
Separable matrices & 1.15, 2.8 \\
Conjugacy classes & 1.16, 2.9 & A070933, A082877
\end{tabular}
\end{center}

\section{The Sequences}

References to the OEIS are accurate as of February 1, 2006. Sequences mentioned in this article may have been added since then and entries in the OEIS may have been modified. 

\subsection{$n \times n$ matrices over $\fq$}
Over the field $\fq$ the number of $n \times n$ matrices is $q^{n^2}$. For $q=2$ this sequence  is A002416 of the OEIS, indexed from $n=0$. The terms for $0 \leq n \leq 4$ are
\[ 1, 2, 16, 512, 65536. \]

\subsection{Invertible matrices}
The number of invertible $n \times n$ matrices is given by 
\[  \gamma_n(q) := |\glnq|=(q^n-1)(q^n - q)(q^n - q^2) \cdots (q^n -q^{n-1}) .\]
We will use $\g_n$ in place of $\g_n(q)$ unless there is a need to be explicit about the base field.
It is convenient to define $\g_0$ to be 1.
 For $q=2$ the sequence  is A002884, and from $\g_0$ to $\g_5$ the terms are
 \[ 1,1,6,168,20160,9999360. \]
A second formula for  $\gamma_n$ is
\[  \gamma_n = (q-1)^n q^{\binom{n}{2}} [n]_q!  ,\]
where ${[n]_q}!:= [n]_q [n-1]_q \cdots [2]_q [1]_q$ and $[i]_q := 1+q+q^2+ \cdots + q^{i-1}$ are the $q$-analogs of $n!$ and $i$.  Also, $[i]_q$ is the number of points in the projective space of dimension $i-1$ over $\fq$. 

A third formula is
\[ \gamma_n= (-1)^n q^{\binom{n}{2}} (q;q)_n  ,\]
where $(a;q)_n$ is defined for $n>0$ by
\[ (a;q)_n := \prod_{j=0}^{n-1}(1-aq^j)  .\]

Construct a random $n \times n$ matrix over $\fq$ by choosing the entries independently and uniformly from $\fq$. Then  $\gamma_n/q^{n^2}$ is the probability that the matrix is invertible, and this probability has a limit as $n \ra \infty$ 
\[ \lim_{n \ra \infty} \frac{\gamma_n}{q^{n^2}} = \prod_{r \geq 1} \left(1-\frac{1}{q^r}\right).\]
For $q=2$ the limit is $0.28878\ldots$. For $q=3$ the limit is $0.56012\ldots$. As $q \ra \infty$ the probability of being invertible goes to 1.

\subsection{Subspaces}
The number of $k$-dimensional subspaces of a vector space of dimension $n$ over $\fq$ is given by 
the  Gaussian binomial coefficient (also called the $q$-binomial coefficient)
\[
\qbinom{n}{k} = \frac{(q^{n} - 1)(q^{n} - q)(q^{n} - q^2) \cdots (q^{n} - q^{k-1}) }
        {(q^{k} - 1)(q^{k} - q)(q^{k} - q^2) \cdots (q^{k} - q^{k-1})}.
\]
Also, we define $\qbinom{0}{0} =1.$ Other useful expressions for the Gaussian binomial coefficients are
\[ \qbinom{n}{k} = \frac{{[n]_q} ! }{{[k]_q}!{[n-k]_q}!} , \]
which shows the $q$-analog nature of the Gaussian coefficients, and
\[ \qbinom{n}{k} = \frac{\g_n}{\g_k \g_{n-k} q^{k(n-k)}} , \]
which comes from the transitive action of $\glnq$ on the set of subspaces of dimension $k$, the denominator being the order of the subgroup stabilizing one of the subspaces.

 Entries A022166--A022188 are the triangles of Gaussian binomial coefficients for $q=2$ to $q=24$. Note that when $q$ is not a prime power the formula does not count subspaces. For $q=2$ the rows from $n=0$ to 6 are
 
\[ \begin{array}{rrrrrrr}
 1&&&&&&\\ 1& 1&&&&& \\ 1& 3& 1& &&&\\1& 7& 7& 1& && \\
1& 15& 35& 15& 1&& \\  1& 31& 155& 155& 31& 1& \\
1& 63& 651& 1395& 651& 63& 1
 \end{array} 
 \]
For the Gaussian binomial coefficients we have the $q$-\textbf{binomial theorem} \cite{Cameron94}
\[  \prod_{1 \leq i \leq n}(1 + q^i t)= \sum_{0 \leq k \leq n}q^{\binom{k}{2}}\qbinom{n}{k} t^k  .\]

Summing over $k$ from 0 to $n$ for a fixed $n$ gives the total number of subspaces. Sequences A006116--A006122 and A015195--A015217  correspond to the values $q=2,\ldots,8$ and $q=9,\ldots,24.$ (Same warning applies to those $q$ that are not prime powers.) For $q=2$ and $n=0,\ldots, 8$ the sequence begins 
\[ 1, 2, 5, 16, 67, 374, 2825, 29212, 417199. \]
\subsection{Splittings}
Let ${n \brace k}_q$ be the number of direct sum splittings of an $n$-dimensional vector space into $k$ non-trivial subspaces without regard to the order among the subspaces. These numbers are $q$-analogs of the Stirling numbers of the second kind, which count the number of partitions of an $n$-set into $k$ non-empty subsets, and the notation follows that of Knuth for the Stirling numbers. Then it is easy to see that
\[ 
   {n \brace k}_{q} = \frac{1}{k!} \sum_{\stackrel{n_{1}+\cdots+n_{k}=n}{n_{i}\geq 1}}
        \frac{\g_{n}}{\g_{n_{1}}\cdots\g_{n_{k}}} .
\]
It is not difficult to verify that the triangle ${n \brace k}_q$ satisfies the two variable generating function identity
\[ 1+ \sum_{n \geq 1}\sum_{k=1}^n {n \brace k}_q \frac{u^n}{\g_n} t^k 
= \exp\left(t \sum_{r \geq 1}\frac{u^r}{\g_r} \right)  .\]
For $q=2$ the entries for $1 \leq n \leq 6$ and $1 \leq k \leq n$ are
\[ \begin{array}{rrrrrr}
   1 &&&&& \\
   1& 3 &&&& \\
   1 & 28 & 28 &&& \\
   1 & 400 & 1680 & 840 && \\  
   1 &10416 & 168640& 277760 & 83328 & \\
   1 & 525792 &  36053248 & 159989760 & 139991040 & 27998208
   \end{array}
 \]

Let $b_n$ be the total number of non-trivial splittings of an $n$-dimensional vector space:
\[ b_n = \sum_{k=1}^n{n \brace k}_q . \]
 The $b_n$ are analogs of the Bell numbers counting the number of partitions of finite sets. Then we have the formula for the generating function for the $b_n$,
\[  1 + \sum_{n \geq 1} b_n \frac{u^n}{\g_n} = \exp\left(\sum_{r \geq 1} \frac{u^r}{\g_r} \right) .\]
For $q=2$ the values of $b_n$ for $n=1,\ldots,6$ are
\[  1, 4, 57, 2921, 540145, 364558049 .\]

Bender and Goldman  \cite{BenderGoldman70} first wrote down the generating function for the $q$-Bell numbers, which, along with the $q$-Stirling numbers, can be put into the context of $q$-exponential families \cite{Morrison04}.

\subsection{Flags}
A \textbf{flag} of length $k$ in a vector space $V$ is an increasing sequence of subspaces
\[  \{0\} = V_0 \subset V_1 \subset V_2 \subset \cdots V_k = V . \] If $\dim V =n$, then a \textbf{complete flag} is a flag of length $n$.
Necessarily, $\dim V_i = i$.

The number of complete flags of an $n$-dimensional vector space over $\fq$ is 
\[ \qbinom{n}{1} \qbinom{n-1}{1} \cdots \qbinom{1}{1} = [n]_q [n-1]_q \cdots [1]_q = [n]_q ! .\]
The reason is that having selected $V_1,\ldots, V_i$, the number of choices for $V_{i+1}$ is $\qbinom{n-i}{1}$. Alternatively, the general linear group acts transitively on the set of complete flags. The stabilizer subgroup of the standard flag is the subgroup of upper triangular matrices, whose order is
 $(q-1)^n q^{\binom{n}{2}}$. Hence, the number of complete flags is
\[ \frac{\g_n}{(q-1)^n q^{\binom{n}{k}}}= [n]_q! .\]

Define $[0]_q ! =1$. The sequence $[n]_q!$ for $q=2$ is A005329. The terms for $n=0,1,2,\ldots,8$ are
\[ 1, 1, 3, 21, 315, 9765, 615195, 78129765, 19923090075. \]

By looking up the phrase ``$q$-factorial numbers" in the OEIS one can find the sequences $[n]_q!$  for $q \leq 13$. The triangle of $q$-factorial numbers is A069777.


The $q$-multinomial coefficient
\[ \qbinom{n}{n_1 \; n_2 \; \cdots \;n_{k}} := \frac{[n]_q!}{[n_1]_q! [n_2]_q! \cdots [n_{k}]_q!}, \]
where $n=n_1 + \cdots + n_{k} $, is the number of flags of length $k$ in an $n$-dimensional space such that $\dim V_i = n_1 + \cdots + n_i$.

\subsection{Linear binary codes}
An $[n,k]$ linear binary code is a $k$-dimensional subspace of the space of $\f2 ^n$. Thus, the number of $[n,k]$ linear binary codes is the Gaussian binomial coefficient $\binom{n}{k}_{\! 2}$ from  section 1.3 and is given by the triangle A022166. 

Two linear binary codes are \textbf{equivalent} (or isometric) if there is a permutation matrix $P$ mapping one subspace to the other. The number of equivalence classes of $[n,k]$ codes is given by the triangle A076831. The early entries are identical with the corresponding binomial coefficients, but that does not hold in general. The rows for $n=0$ to 6 are
\[ \begin{array}{rrrrrrr}
 1&&&&&&\\ 1& 1&&&&& \\ 1& 2& 1& &&&\\1& 3& 3& 1& && \\
1& 4& 6& 4& 1&& \\  1& 5& 10& 10& 5& 1& \\
1& 6& 16& 22& 16& 6& 1
 \end{array} 
 \]
For linear codes over other finite fields, the notion of equivalence uses matrices $P$ having exactly one non-zero entry in each row and column. These matrices, as linear maps,  preserve the Hamming distance between vectors. They form a group isomorphic to the wreath product of the multiplicative group of $\fq$ with $S_n$. 

There are more than a dozen sequences and triangles associated to linear binary codes in the OEIS. They are listed in the short index of the OEIS under ``Codes.''

\subsection{Matrices by rank}
The number of $m \times n$ matrices of rank $k$ over $\fq$ is
 \begin{align*}
	&  \qbinom{m}{k} (q^{n}-1)(q^{n}-q)\cdots(q^{n}-q^{k-1})\\
	&=\qbinom{n}{k}(q^{m}-1)(q^{m}-q)\cdots(q^{m}-q^{k-1}) \\
	&= \frac{(q^{m}-1)(q^{m}-q)\cdots(q^{m}-q^{k-1})\,\,\,(q^{n}-1)
	            (q^{n}-q)\cdots(q^{n}-q^{k-1})}
		    {(q^{k} - 1)(q^{k} - q) \cdots (q^{k} - q^{k-1})}.
    \end{align*}
(To justify the first line of the formula note that the number of $k$-dimensional subspaces of $\fq^m$ to serve as the column space of a rank $k$ matrix is $\qbinom{m}{k}$. Identify the column space with the image of the associated linear map from $\fq^n$ to $\fq^m$. There are $(q^{n}-1)(q^{n}-q)\cdots(q^{n}-q^{k-1})$ surjective linear maps from $\fq^n$ to that $k$-dimensional image. The second line follows by transposing.)

Define the triangle  $r(n,k)$ to be the number of $n \times n$ matrices of rank $k$ over $\fq$ and  $r(0,0) := 1$.  Thus,
\[  r(n,k) = \frac{\left((q^{n}-1)
	            (q^{n}-q)\cdots(q^{n}-q^{k-1})\right)^2}
		    {(q^{k} - 1)(q^{k} - q) \cdots (q^{k} - q^{k-1})}.
\]
For $q=2$ the entries from $r(0,0)$ to $r(5,5)$ are
\[ \begin{array}{rrrrrr}
   1 &&&&& \\
   1& 1 &&&& \\
   1 & 9 & 6 &&& \\
   1 & 49 & 294 & 168 && \\
   1 & 225 & 7350 & 37800 & 20160 \\
   1 & 961 & 144150 & 4036200 & 19373760 & 9999360 
   \end{array}
 \]

\subsection{Linear and projective derangements}
A matrix is a \textbf{linear derangement} if it is invertible and does not fix any non-zero vector. Such a matrix is characterized as not having 0 or 1 as an eigenvalue.   Let $e_n$ be the number of linear derangements and define $e_0 = 1$. Then $e_n$ satisfies the recursion
\[ e_n = e_{n-1}(q^n-1)q^{n-1} + (-1)^n q^{n(n-1)/2} .\]
For $q=2$, the sequence is A002820 (with offset 2) and the first few terms beginning with $e_0$ are  $1,0, 2, 48, 5824, 2887680, \ldots$. The sequence can be obtained from the generating function
\[   1 + \sum_{n \geq 1} \frac{e_n}{\gamma_n} u^n = 
         \frac{1}{1-u}\prod_{r \geq 1}\left( 1- \frac{u}{q^r} \right) .\]
The proof of this is in section 2. The asymptotic probability that an invertible matrix is a linear derangement is 
\[   \lim_{n \ra \infty} \frac{e_n} {\gamma_n} = \prod_{r \geq 1} \left(1-\frac{1}{q^r}\right).\]
The asymptotic probability that a square matrix is a linear derangement is
\[ \lim_{n \ra \infty} \frac{e_n} {q ^{n^2}} = \prod_{r \geq 1} \left(1-\frac{1}{q^r}\right)^2
 \]

A matrix over $\fq$ is a \textbf{projective derangement} if the induced map on projective space has no fixed points. Equivalently, this means that the matrix has no eigenvalues in $\fq$.   If $d_n$ is the number of such $ n \times n$ matrices, then 
\[  1 + \sum_{n \geq 1} \frac{d_n}{\gamma_n} u^n = 
         \frac{1}{1-u}\prod_{r \geq 1}\left( 1- \frac{u}{q^r} \right)^{q-1}  .\]
For $q=2$ the notions of projective and linear derangement are the same, and the corresponding sequence is given above. For $q=3$ the sequence 
has initial terms from $n=1$ to 5
\[0, 18, 3456, 7619508, 149200289280.\]
Since two matrices that differ by a scalar multiple induce the same map on projective space, the number of maps that are projective derangements is $d_n/(q-1)$. The asymptotic probability that a random invertible matrix is a projective derangement is the limit 
\[  \lim_{n \ra \infty} \frac{d_n} {\gamma_n} = \prod_{r \geq 1} \left(1-\frac{1}{q^r}\right)^{q-1}.\]
The asymptotic probability that a random $n \times n$ matrix is a projective derangement is
\[  \lim_{n \ra \infty} \frac{d_n} {q ^{n^2}} = \prod_{r \geq 1} \left(1-\frac{1}{q^r}\right)^q.\]
As $q \ra \infty$, this probability goes to $1/e$, the same value as the asymptotic probability that a random permutation is a derangement.

\subsection{Diagonalizable matrices}
In this section let $d_n$ be the number of diagonalizable $n \times n$ matrices over $\fq$. Then
\[ 1 + \sum_{n \geq 1} \frac{d_n}{\g_n}u^n = \left(\sum_{m \geq 0}\frac{u^m}{\g_m} \right)^q .\]
It follows that 
\[  d_n = \sum_{n_1 + \cdots + n_q = n} \frac{\g_n}{\g_{n_1}\cdots\g_{n_q}} . \]
 For $q=2$ this simplifies to
 \[  d_n = \sum_{i=0}^n \frac{\g_n}{\g_i \g_{n-i}}  .\]
 The sequence for  $d_1$ to $d_8$ is
 \[ 2, 8, 58, 802, 20834, 1051586, 102233986, 196144424834.\]
 For $q=3$ the five initial terms are 
 \[ 3, 39, 2109, 417153, 346720179. \]  
 
 For arbitrary $q$, we can easily find $d_2$,
 \begin{align*} d_2  &= \sum_{n_1 + \cdots + n_q = 2} \frac{\g_n}{\g_{n_1}\cdots\g_{n_q}}  \\
                                           &= q \frac{\g_2}{\g_0 \g_2} + {\binom{q}{2}} \frac{\g_2}{\g_1 \g_1} \\
                                           &= \frac{q^4-q^2+2q}{2} .
 \end{align*}

\subsection{Projections}
A \textbf{projection} is a matrix $P$ such that $P^2=P$. Let $p_n$ be the number of $n \times n$ projections. Then
\[
 1+\sum_{n\geq 1} \frac{p_n}{\g_n}u^n= \left( \sum_{m \geq 0} \frac{u^m}{\g_m} \right)^2.
\]
It follows that 
\[  p_n = \sum_{0 \leq i \leq n} \frac{\g_n}{\g_i \g_{n-i}} . \]
In the sum the term $\g_n / \g_i \g_{n-i} $ is the number of projections of rank $i$.
Projections are also characterized as  diagonalizable matrices having eigenvalues 0 or 1. Thus, for $q=2$ the diagonalizable matrices are precisely the projections, and so we get the same sequence as in section 1.9. From $p_1$ to $p_8$ the sequence is
\[  2, 8, 58, 802, 20834, 1051586, 102233986, 196144424834.\]

For $q=3$ the map $P \mapsto P + I$ gives a bijection from the set of projections to the set of diagonalizable matrices with eigenvalues 1 or 2. Such matrices are precisely the solutions of $X^2=I$. The sequence is given in the OEIS by  A053846. From $p_1$ to $p_7$ the sequence is
\[1, 2, 14, 236, 12692, 1783784, 811523288, 995733306992.\]

There is a bijection between the set of $n \times n$ projections and the direct sum splittings $\fq^n = V \oplus W$. The projection $P$ corresponds to the splitting $\textrm{Im} P \oplus \textrm{Ker} P$. Note that $V \oplus W$ is regarded as different from $W \oplus V$. Hence, 
\[  p_n = 2 + 2 {n \brace 2}_q ,\]
because ${n \brace 2}_q$ is the number of splittings into two proper subspaces, with $V \oplus W$ regarded as the same as $W \oplus V$.

From \S 1.3 we have
\[ \qbinom{n}{k} = \frac{\g_n}{\g_k \g_{n-k} q^{k(n-k)}} , \]
showing the relationship between the number of subspaces of dimension $k$ and the number of projections of rank $k$. From it we see that t there are $q^{k(n-k)}$ complementary subspaces for a fixed subspace of dimension $k$.

\subsection{Solutions of $A^k=I$ }

Let $a_n$ be the number of $n \times n$ matrices $A$ satisfying $A^2=I$. 
Such matrices  correspond to group homomorphisms from the cyclic group of order $2$ to $\glnq$. 
 In characteristic other than two, the generating function for the $a_n$ is
\[  1+ \sum_{n \geq 1} \frac{a_n}{\g_n}u^n = \left( \sum_{m \geq 0}\frac{u_m}{\g_m} \right)^2, \]
and so 
\[  a_n = \sum_{i=0}^{n} \frac{\g_n}{\g_i \g_{n-i} }. \]
These matrices are those which are diagonalizable and have only 1 and $-1$ for eigenvalues. If we fix two distinct elements $\lambda_1$ and $\lambda_2$ in the base field $\fq$, then $a_n$ also gives the number of diagonalizable matrices having eigenvalues $\lambda_1$ and $\lambda_2$. Taking the eigenvalues to be 0 and 1 gives the set of projections, and so we see that the number of projections is the same as the number of solutions to $A^2=I$, when $q$ is not a power of 2. The proof in \S 2.4 for the number of projections is essentially the proof for the more general case. The sequence A053846 gives the number of solutions over $\F_3$.

In characteristic two, $A^2=I$ does not imply that $A$  is diagonalizable and the previous formula does not hold. 
We have the formula
\[  a_n =  \sum_{0 \leq i \leq n/2} \frac{\g_n}{q^{i(2n-3i)}\g_{i} \g_{n-2i} }\]
to be proved in \S 2.5. Because $A^2=I$ is equivalent to $(A+I)^2=0$, the formula for $a_n$ also counts the number of $n \times n$ nilpotent matrices $N$ such that $N^2=0$.
For $q=2$ the sequence is A053722. Although there appears to be an error in the formula given in the OEIS entry, the initial terms of the sequence are correct. For  $n=1,\ldots,8$ they are
\[ 1,4,22,316,6976,373024,32252032,6619979776. \]
For $q=4$ the sequence is A053856 and the initial terms for $n=1,\ldots,7$ are
\[1,16,316,69616,21999616,74351051776,374910580965376.\]

More generally, let $k$ be a positive integer not divisible by $p$, where $q$ is a power of $p$. We consider the solutions of $A^k=I$ and let $a_n$ be the number of $n \times n$ solutions with coefficients in $\fq$. Now $z^k-1$ factors into a product of distinct irreducible polynomials
\[ z^k-1 = \phi_1(z)\phi_2(z) \cdots \phi_r(z) .\]
Let $d_i = \deg \phi_i$. 
Then the generating function for the $a_n$ is 
\[ 1+ \sum_{n \geq 1} \frac{a_n}{\g_n} u^n = \prod_{i=1}^r \sum_{m \geq 0} \frac{q^{md_i}}{\g_m(q^{d_i})} .\]
Note that $a_n$ counts the homormorphisms from the cyclic group of order $k$ to $\glnq$. The generating function given here is a special case of the generating function given by Chigira, Takegahara, and Yoshida \cite{CTY00} for the sequence whose $n$th term is the number of homomorphisms from a finite group $G$ to $\glnq$ under the assumption that the characteristic of $\fq$ does not divide the order of $G$.

Now we consider some specific examples. Let $k=3$ and let $q$ be a power of 2. Then \[z^3-1=z^3+1=(z+1)(z^2+z+1) \] is the irreducible factorization. Thus, $d_1=1$ and $d_2=2$. The generating function is
\[   1+ \sum_{n \geq 1} \frac{a_n}{\g_n}u^n= \left(\sum_{m \geq 0}\frac{u^m}{\g_m}  \right)
                    \left(  \sum_{m \geq 0}\frac{u^{2m}}{\g_m(q^2)}  \right).
\]                    
With $q=2$ we get the sequence A053725 with initial terms for $n=1,\ldots,8$
\[ 1,3,57,1233,75393,19109889,6326835201,6388287561729.  \]
With $q=4$ we get the sequence A053857 with initial terms for $n=1,\ldots, 5$
\[  3,63,8739,5790339,25502129667  .\]
For the next example, let $k=8$ and let $q$ be a power of 3. Then \[z^8-1=(z-1)(z+1)\phi_3(z) \phi_4(z) \phi_5(z), \] where the last three factors all have degree 2. Thus, $r=5$, $d_1=d_2=1$, and 
$d_3=d_4=d_5=2$. The generating function is given by
\[    \sum_{n \geq 1} \frac{a_n}{\g_n}u^n= \left(\sum_{m \geq 0}\frac{u^m}{\g_m}  \right)^2
                    \left(  \sum_{m \geq 0}\frac{u^{2m}}{\g_m(q^2)}  \right)^3.
\]
For $q=3$ this gives the sequence A053853 beginning for $n=1,\ldots,5$
\[ 2,32,4448,3816128,26288771456  .\]
These are the sequences in the OEIS for small values of $k$ and $q$.
\begin{center}
\begin{tabular}{rccc}
$k$ & $q=2$ &  $q=3$ & $q=4$ \\
2 &  A053722 & A053846 &  A053856\\
 3 & A053725  & A053847 & A053857\\
 4  & A053718 & A053848 & A053859\\
 5 &  A053770 & A053849 & A053860\\ 
 6 & A053771  & A053851 & A053861\\ 
 7 & A053772  & A053852 & A053862\\ 
 8 &  A053773 & A053853 & A053863\\ 
 9 & A053774  & A053854 &\\ 
10 & A053775 & A053855 &\\ 
11 &  A053776 & & \\ 
12 &  A053777 & & \\
\end{tabular}
\end{center}

\subsection{Nilpotent matrices}
Fine and Herstein \cite{FineHerstein58} proved that the number of nilpotent $n \times n$ matrices is $q^{n(n-1)}$, and Gerstenhaber \cite{Gerstenhaber61} simplified the proof.  Recently Crabb has given an even more accessible proof \cite{Crabb06}. For $q=2$ this is sequence  A053763 with initial terms from $n=1$ to $n=6$
\[1,1,4,64,4096,1048576,1073741824. \] 

\subsection{Cyclic matrices} 
A matrix $A$ is \textbf{cyclic} if there exists a vector $v$ such that $\{ A^i v | i=0,1,2,\ldots \}$ spans the underlying vector space. (The term \textbf{regular} is also used. ) An equivalent description is that the minimal and characteristic polynomials of $A$ are the same. Let $a_n$ be the number of cyclic matrices over $\fq$. The generating function factors as
\[ 1+ \sum_{n \geq 1} \frac{a_n}{\g_n}u^n = \prod_{d \geq 1}
       \left(1 + \frac{1}{q^d-1} \frac{u^d}{1-(u/q)^d} \right) ^{\nu_d} ,\] 
and can be put into the form
\[ 1+ \sum_{n \geq 1} \frac{a_n}{\g_n}u^n = \frac{1}{1-u}\prod_{d \geq 1}
       \left(1 + \frac{u^d}{q^d(q^d-1)} \right) ^{\nu_d} ,\] 
where $\nu_d$ is the number of irreducible, monic polynomials of degree $d$ over $\fq$. The generating function can be extracted from the proof of Theorem 1 in \cite{Fulman02}, which we present in section 2. 

For $q=2$ the sequence from $a_1$ to $a_7$  is 
\[ 2, 14, 412, 50832, 25517184, 51759986688, 422000664182784 .\]
Wall \cite{Wall99} proved that the probability that an $n \times n$ matrix is cyclic has a limit 
\[  \lim_{n \ra \infty} \frac{a_n}{q^{n^2}} = \left(1 - \frac{1}{q^5} \right) \prod_{r \geq 3} 
                   \left( 1-  \frac{1}{q^r} \right)  .
\] Fulman \cite{Fulman02} also has a proof.
For $q=2$ this limit is $0.7403\ldots$.    

\subsection{Semi-simple matrices}
A matrix $A$ is \textbf{semi-simple} if it diagonalizes over the algebraic closure of the base field. Let $a_n$ be the number of semi-simple $n \times n$ matrices over $\fq$. Then the generating function has a factorization
\[ 1 + \sum_{n \geq 1} \frac{a_n}{\g_n}u^n =  
          \prod_{d \geq 1}\left( 1 + \sum_{j \geq 1} \frac{u^{jd}}{\g_j(q^d) }\right)^{\nu_d},
\]
where $\g_j(q^d) = |\gl_j(q^d) |.$
For $q=2$ the sequence from $a_1$ to $a_7$ is
\[   2, 10, 218, 25426, 11979362, 24071588290, 195647202043778.\]

\subsection{Separable matrices}     
A matrix is \textbf{separable} if it is both cyclic and semi-simple, which is equivalent to having a characteristic polynomial that is square-free. Let $a_n$ be the number of separable $n \times n$ matrices over $\fq$. Then the generating function factors
\[ 
    1 + \sum_{n \geq 1} \frac{a_n}{\g_n} u^n =
    \prod_{d \geq 1} \left( 1 + \frac{u^d}{q^d-1} \right)^{\nu_d}.
\]     
This can also be factored as
\[  
 1 + \sum_{n \geq 1} \frac{a_n}{\g_n} u^n =  \frac{1}{1-u}
    \prod_{d \geq 1} \left( 1 + \frac{u^d(u^d -1)}{q^d(q^d-1)} \right)^{\nu_d}.
\]     
 
For $q=2$ the sequence from $a_1$ to $a_7$ is 
\[
2,  8,  160,  22272,  9744384,  20309999616,  165823024988160.
\]

The number of conjugacy classes of separable $n \times n$ matrices is $q^n -q^{n-1}$ for $n \geq 2$ and $q$ for $n=1$. This is proved by Neumann and Praeger \cite[Lemma 3.2]{NeumannPraeger95} by showing that the number of square-free monic polynomials of degree $n$ is $q^n - q^{n-1}$. There is a natural bijection between the set of conjugacy classes of separable matrices and the set of square-free monic polynomials obtained by associating to each conjugacy class the characteristic polynomial of the matrices in that class.

\subsection{Conjugacy classes}
Let $a_n$ be the number of conjugacy classes of $n \times n$ matrices over $\fq$. The generating function for this sequence is
\[   1 + \sum_{n \geq 1} a_n u^n = \prod_{r \geq 1} \frac{1}{1-qu^r}   .\]
The number of conjugacy classes grows like $q^n$. In fact,
\[  \lim_{n \ra \infty} \frac{a_n}{q^n} = \prod_{r \geq 1}\left( 1- \frac{1}{q^r} \right)^{-1} , \]
which is the reciprocal of the limiting probability that a matrix is invertible. See section 1.2.

For $q=2$ the sequence is A070933.
The initial terms for $n=1,\ldots,10$ are 
\[2, 6, 14, 34, 74, 166, 350, 746, 1546, 3206. \]

For $q=3$ the initial terms for $n=1,\ldots,10$ are
\[3, 12, 39, 129, 399, 1245, 3783, 11514, 34734, 104754. \]

Let $b_n$ be the number of conjugacy classes in the general linear group $\glnq$, the group of $n \times n$ invertible matrices over $\fq$. Then
\[ 1+ \sum_{n \geq 0} b_n u^n = \prod_{r \geq 1} \frac{1-u^r}{1-q u^r}  .\]
 
 For $q=2$ the sequence is  A006951, which starts  \[ 1, 3, 6, 14, 27, 60, 117, 246, 490, 1002. \]
 
 For $q=3$ the sequence is A006952, which starts \[ 2, 8, 24, 78, 232, 720, 2152, 6528, 19578, 58944. \]
 
 For $q=4, 5, 7$ the sequences are A049314, A049315, and A049316.
 
 The number of conjugacy classes is asymptotic to $q^n$:
 \[ \lim_{n \ra \infty} \frac{b_n}{q^n} =  1. \]
 Hence, in the limit the ratio of the number of conjugacy classes of invertible matrices to the number of conjugacy classes of all matrices is the same as the limiting probability that a matrix is invertible. That is, \[ \lim_{n \ra \infty} \frac{b_n}{a_n} = \prod_{r \geq 1} \left(1- \frac{1}{q^r} \right) .\]
 
Sequence A070731 gives the size of the largest conjugacy class in $\gl_n(2)$. Starting with $n=1$ the initial terms are  \[1, 3, 56, 3360, 833280, 959938560.\] The minimal order of the centralizers of elements in $\gl_n(2)$ is given by the quotient of $\g_n$ by the $n$th term in this sequence. The resulting sequence for $n=1,\ldots 10$ is 
\[1,  2,  3,  6,  12,  21,  42,  84, 147, 294. \]
This sequence is A082877 in the OEIS. 

\section{Selected Proofs}

\subsection{The cycle index and generating functions}
In sections 1.6-1.13 we make heavy use of generating functions of the form
\[  1 + \sum_{n \geq 1} \frac{a_n}{\g_n}u^n ,\]
where the sequence $a_n$ counts some class of $n \times n $ matrices. These generating functions come from the \textbf{cycle index} for matrices that was first defined by Kung \cite{Kung81} and later extended by Stong \cite{Stong88}. We follow Fulman's notation in \cite{Fulman02}. The cycle index for conjugation action of $\glnq$ on $\mnq$ is a polynomial in the  indeterminates  $x_{\phi, \lambda}$, where $\phi$ ranges over the  set $\Phi$ of monic irreducible polynomials over $\fq$ and $\lambda$ ranges over the partitions of the positive integers. First we recall that the conjugacy class of a matrix $A$ is determined by the isomorphism type of the associated $\fq[z]$-module on the vector space $\fq^n$ in which the action of $z$ is that of $A$. This module is isomorphic to a direct sum 
\[ \bigoplus_{i=1}^{k}\bigoplus_{j=1}^{l_{i}} 
   \fq[z] / (\phi_{i}^{\lambda_{i,j}})
\]
where $\phi_{1},\dots,\phi_{k}$ are distinct monic irreducible 
polynomials; for each $i$, $\lambda_{i,1} \geq \lambda_{i,2}\geq 
\dots \geq \lambda_{i,l_{i}}$ is a partition of $n_{i} = 
\sum_{j}\lambda_{i,j}$. Since $A$ is $n \times n$, then 
$n=\sum_{i}n_{i} \deg \phi_{i} = \sum_{i,j}\lambda_{i,j}\deg \phi_{i}$.
Let $\lambda_{i}$ denote the partition of $n_{i}$ given by the 
$\lambda_{i,j}$ and define $|\lambda_{i}|=n_{i}$. The conjugacy class 
of $A$ in $\mnq$ is determined by the data consisting of the 
finite list of distinct monic irreducible polynomials 
$\phi_{1},\dots,\phi_{k}$ and the partitions 
$\lambda_{1},\dots,\lambda_{k}$.

The cycle index is defined to be
\[   \frac{1}{\g_n} \sum_{A \in \mnq} \prod_{\phi \in \Phi} x_{\phi, \lambda_\phi (A)} ,\]
where $ \lambda_\phi (A)$ is the partition associated to $\phi$ in the conjugacy class data for $A$. If $\phi$ does not occur in the polynomials associated to $A$, then $\lambda_\phi(A)$ is the empty partition, and we define $x_{\phi, \lambda_\phi (A)} =1 $.

We construct the generating function for the cycle index 
\[ 1 + \sum_{n \geq 1} \frac{u^n}{\g_n} \sum_{A \in \mnq} \prod_{\phi \in \Phi} x_{\phi, \lambda_\phi (A)}.
 \]
This generating function has a formal factorization over the monic irreducible polynomials. For the proof we refer to the paper of Stong \cite[p. 170]{Stong88} .
\begin{lemma}
 \[ 1 + \sum_{n=1}^{\infty}\frac{u^{n}}{\g_n}\sum_{A \in \mnq}
   \prod_{\phi}x_{\phi,\lambda_{\phi}(A)}  
   = \prod_{\phi}\sum_{\lambda}\frac{x_{\phi,\lambda}u^{|\lambda|\deg 
   \phi}}{c_{\phi}(\lambda)}
   \]
   where $c_{\phi}(\lambda)$ is the order of the group of module 
   automorphisms of the $\fq[z]$-module 
   $\bigoplus_{j}\fq[z]/(\phi^{\lambda_{j}})$.
\end{lemma}

\begin{lemma}\label{lemma.b}
  Fix a monic irreducible $\phi$. Then
  \[ \sum_{\lambda}\frac{u^{|\lambda| \deg \phi}}{\cpl} = 
     \prod_{r \geq 1}\left(1-\frac{u^{\deg \phi}}{q^{r \deg 
     \phi}}\right) ^{-1}.
  \]
\end{lemma}
\proof We use the formula for $\cpl$ proved by Kung \cite[Lemma 2, p. 146]{Kung81}. Let $\lambda$ be a partition of $n$ and let $b_{i}$ be the number of parts of size 
 $i$. Let $d_{i} = b_{1} +2b_{2}+ \dots + ib_{i} +i(b_{i+1}+\dots + b_{n})$.
 Then
 \[ c_{\phi}(\lambda) = 
 \prod_{i}\prod_{k=1}^{b_{i}}(q^{d_{i}\deg 
 \phi}-q^{(d_{i}-k)\deg \phi}).
 \]  We see from this formula that $\cpl$ is a  function of $\lambda$ and $q^{\deg \phi}$. Therefore it is sufficient to prove the lemma for a single polynomial of degree 1, say 
$\phi(z)=z$, and then to replace $u$ by $u^{\deg \phi}$ and $q$ by $q^{\deg \phi}$.  So, now we will prove that for $\phi(z)=z$,
 \[ \sum_{\lambda}\frac{u^{|\lambda|}}{\cpl} = 
     \prod_{r \geq 1}\left(1-\frac{u}{q^{r }}\right) ^{-1}.
  \]
We split the left side into an outer sum over $n$ and an inner sum over 
the partitions of $n$
\[ 
 \sum_{\lambda}\frac{u^{|\lambda|}}{\cpl}  = 
 1 + \sum_{n \geq 1} \sum_{|\lambda|=n }\frac{u^{n}}{\cpl}.
\] 
To evaluate the inner sum we note that $\g_n\sum_{|\lambda|=n}\frac{1}{\cpl}$ 
is the number of $n \times n$ nilpotent matrices, which is $q^{n(n-1)}$ from the theorem of Fine and 
Herstein \cite{FineHerstein58}.
Therefore, the $u^n$ coefficient of the left side is
\[
 \sum_{|\lambda|=n}\frac{1}{\cpl}  =
 \frac{1}{q^{n}(1-\frac{1}{q})\dots(1-\frac{1}{q^{n}})}.
\]
For the right side we use a formula of Euler, which is a special case of Cauchy's identity and a limiting case of the $q$-binomial theorem. A convenient reference is the book of Hardy and Wright 
\cite[Theorem 349, p. 280]{HardyWright79}. For $|a|<1$, $|y|<1$, the 
coefficient of $a^{n}$ in the infinite product 
\[\prod_{r \geq 1}\frac{1}{(1-ay^{r})} \]
is 
\[ \frac{y^{n}}{(1-y)(1-y^{2})\dots(1-y^{n})} .\]
Let $a=u$ and $y=\frac{1}{q}$ to see that the 
$u^{n}$ coefficient of
\[ \prod_{r \geq 1}\left( 1-\frac{u}{q^{r}} \right)^{-1} \]
is \[  \frac{1}{q^{n}(1-\frac{1}{q})\dots(1-\frac{1}{q^{n}})}.\] 
Therefore, the coefficients of $u^n$ of the left and right sides are equal and so we have proved that for $\phi(z)=z$,
 \[ \sum_{\lambda}\frac{u^{|\lambda|}}{\cpl} = 
     \prod_{r \geq 1}\left(1-\frac{u}{q^{r }}\right) ^{-1}.
  \] \qed

\begin{lemma}
 Let $\Phi ' $ be a subset of the irreducible monic polynomials. Let $a_n$ be the number of $n \times n$ matrices whose conjugacy class data involves only polynomials $\phi \in \Phi'$. Then
 \begin{align*}  1 + \sum_{n \geq 1}\frac{a_n}{\g_n} u^n &=  \prod_{\phi \in \Phi'}
      \sum_{\lambda}\frac{u^{|\lambda| \deg \phi}}{\cpl} \\
 1 + \sum_{n \geq 1}\frac{a_n}{\g_n} u^n & =     \prod_{\phi \in \Phi'} \prod_{r \geq 1}
               \left(1- \frac{u^{\deg \phi}}{q^{r \deg \phi}} \right) ^{-1}.
\end{align*}           
 \end{lemma}
\proof
Using Lemma 1 we set $x_{\phi,\lambda}=1$ or $0$ according to whether $\phi \in \Phi'$ or not. Then on the left side the inner sum is over the matrices $A$ that do not have factors of $\phi \in \Phi'$ in their characteristic polynomials, and so the coefficient of $u^n$ is simply the number of such $n \times n$ matrices. This gives the first equality in the statement of the lemma
\[ 1 + \sum_{n \geq 1}\frac{a_n}{\g_n} u^n = \prod_{\phi \in \Phi'}
      \sum_{\lambda}\frac{u^{|\lambda| \deg \phi}}{\cpl}  .
\]
Then from Lemma 2 we get the second equality.  \qed      
 
\begin{lemma}
\begin{align*}  
   \frac{1}{1-u}&=\prod_{\phi \neq z}\sum_{\lambda}\frac{u^{|\lambda | \deg \phi}}{\cpl} \\
   \frac{1}{1-u}&=\prod_{\phi \neq z}\prod_{r \geq 1}
               \left(1- \frac{u^{\deg \phi}}{q^{r \deg \phi}} \right) ^{-1}.
\end{align*}               
\end{lemma}
\proof Using Lemma 3 we let $\Phi'$ be the complement of the polynomial $\phi(z)=z$. The matrices not having factors of $z$ in their characteristic polynomials are the invertible matraices. Thus,  $a_n = \g_n$, and so the left side is $1 + \sum_{n \geq 1} u^n$. 
\qed

We can generalize the first part of Lemma 3 to allow for conjugacy class data in which the allowed partitions vary with the polynomials. The proof follows immediately from Lemma 1.
\begin{lemma}
For each monic, irreducible polynomial $\phi$ let $L_\phi$ be a subset of all partitions  of the positive integers. Let $a_n$ be the number of $n \times n$ matrices $A$ such that $\lambda_\phi(A) \in L_\phi$ for all $\phi$. Then
\[  1+\sum_{n \geq 1} \frac{a_n}{\g_n}u^n = \prod_{\phi}\sum_{\lambda \in L_\phi} \frac{u^{|\lambda| \deg \phi}}{\cpl} .
\]
\end{lemma}

\begin{lemma}
\[
  {1-u} = \prod_{\phi} \left(1-\frac{u^{\deg \phi}} {q^{\deg \phi} }\right)=
  \prod_{d \geq 1} \left(1-\frac{u^d}{q^d} \right)^{\nu_d}
\]
\end{lemma}
\proof 
Unique factorization in the ring of polynomials $\fq[z]$ means that each of the $q^n$ monic polynomials of degree $n$ has a unique factorization as a product of monic irreducible polynomials. This implies the factorization of the generating function
\[  \sum_{n \geq 0} q^n u^n = \prod_{\phi}\sum_{k \geq 0} u^{k \deg \phi} .\] 
The left side and the inner sum on the right are geometric series, and so
\[ \frac{1}{1-qu} = \prod_{\phi}  \frac{1}{1-u^{\deg \phi}}. \] 
Grouping the factors on the right according to degree we get
\[ \frac{1}{1-qu} = \prod_{d \geq 1} \left( \frac{1}{1-u^d} \right)^{\nu_d}. \] 
Substituing $u/q$ for $u$ gives
\[ \frac{1}{1-u} = \prod_{d \geq 1} \left( \frac{1}{1-(u/q)^d} \right)^{\nu_d}. \] 
The lemma follows by taking reciprocals.
\qed

The final lemma below is used to evaluate limiting probabilities such as $\lim_{n \ra \infty} {a_n/\g_n}$. The proof is straightforward and omitted.
\begin{lemma}
  If \[ \sum_{n \geq 0} \alpha_n u^n = \frac{1}{1-u}F(u) \]
  where $F(u)$ is analytic and the series for $F(1)$ is convergent, then
  \[  \lim_{n \ra \infty} \alpha_n =F(1)  .\]
\end{lemma}
Now with these lemmas we are ready to prove several results stated in \S 1.

\subsection{Linear and Projective Derangements}
In this section we prove the results given in \S 1.8.  Linear derangements are matrices with no eigenvalues of 0 or 1, which means that their characteristic polynomials do not have factors of $z$ or $z-1$. 
\begin{theorem}
Let $e_n$ be the number of $n \times n$ linear derangements. Then
\[ 1 + \sum_{n \geq 1} \frac{e_n}{\g_n}u^n = \frac{1}{1-u} \prod_{r \geq 1}
      \left(1-\frac{u}{q^r} \right) .
\]      
\end{theorem}
\proof
In Lemma 3 let $\Phi'=\Phi \setminus \{z,z-1\}$ to see that
\[ 1 + \sum_{n \geq 1} \frac{e_n}{\g_n}u^n = \prod_{\phi \in \Phi'} \prod_{r \geq 1}
               \left(1- \frac{u^{\deg \phi}}{q^{r \deg \phi}} \right) ^{-1}.
\]
On the right side multiply and divide by the product corresponding to $\phi(z)=z$, which is 
 \[  \prod_{r \geq 1} \left(1- \frac{u}{q^r} \right)^{-1} , \]
 to see that
\[1 + \sum_{n \geq 1} \frac{e_n}{\g_n}u^n = \prod_{\phi \neq z} \prod_{r \geq 1}
               \left(1- \frac{u^{\deg \phi}}{q^{r \deg \phi}} \right) ^{-1} \, 
               \prod_{r \geq 1} \left(1- \frac{u}{q^r} \right).
\]
Then use Lemma 4 to give
\[ 1 + \sum_{n \geq 1} \frac{e_n}{\g_n}u^n = \frac{1}{1-u}\, 
               \prod_{r \geq 1} \left(1- \frac{u}{q^r} \right).
\]
\qed
 
From the generating function we can derive a recursive formula for $e_n$.
\begin{corollary}
  \[ e_n = e_{n-1}(q^n-1)q^{n-1} + (-1)^n q^{n(n-1)/2}, \, e_0 = 1.\]
\end{corollary}
\proof
 From Theorem 8 it follows that $e_n/\g_n$ is the sum of the $u^i$ coefficients of 
 \[ \prod_{r \geq 1}(1 - u/q^r)\]  for $i=0,1,\ldots,n$. Now the $u^i$ coefficient is
$$ (-1)^i \sum_{ 1 \leq r_1 < r_2 < \cdots < r_i} \frac{1}{q^{r_1 + r_2 + \cdots + r_i}}. $$ By induction one can easily show that this coefficient is
$$  \frac{(-1)^i }{(q^i -1)(q^{i-1}-1) \cdots (q-1)}.$$
Therefore
$$  \frac{e_n}{\g_n} = 1 + \sum_{1 \leq i \leq n}  \frac{(-1)^i }{(q^i -1)(q^{i-1}-1) \cdots (q-1)}. $$
Next,
$$  \frac{e_n}{\g_n} = \frac{e_{n-1}}{\g_{n-1}} + \frac{(-1)^n}{(q^n -1)\cdots (q-1)} .$$
Making use of the formula for $\g_n$ and $\g_{n-1}$ and canceling where possible we see that
$$e_n = e_{n-1}(q^n-1)q^{n-1} + (-1)^n q^{n(n-1)/2}. $$
\qed
 
We present the next result because the sequence $e_n$ for $q=2$ is given in the OEIS as $2^{n(n-1)/2 }a_n$ where $a_n$ satisfies a recursive formula.
\begin{corollary}
  Let 
  \[ a_n = \frac{e_n}{q^{n(n-1)/2}}. \]
  Then $a_n$ satisfies the recursion
  \[  a_n=a_{n-1}(q^n-1)+(-1)^n, \, a_0=1.\]
\end{corollary}
\proof  The proof follows immediately from the recursive formula for $e_n$.
\qed

\begin{corollary} The asymptotic probability that an invertible matrix is a linear derangement is
  \[ \lim_{n \ra \infty} \frac{e_n}{\g_n} = \prod_{r \geq 1}\left(1 - \frac{1}{q^r} \right) . \]
The asymptotic that any matrix is a linear derangement is 
  \[ \lim_{n \ra \infty} \frac{e_n}{q^{n^2}} =  \prod_{r \geq 1}\left(1 - \frac{1}{q^r} \right)^2 .\]
\end{corollary}
\proof
We use Lemma 7 for the first statement. Then
 \begin{align*}
   \lim_{n \ra \infty} \frac{e_n}{q^{n^2}}&=
       \lim_{n \ra \infty} \frac{e_n}{\g_n} \frac{\g_n}{q^{n^2}}  \\
     &=\left( \lim_{n \ra \infty}\frac{e_n}{\g_n}\right) \left( \lim_{n \ra \infty} \frac{\g_n}{q^{n^2}}\right) . 
\end{align*}
We have just computed the first limit on the right and the second limit (from \S 1.2) is the same.
\qed

The proof of Theorem 8 easily adapts to give the generating function for the number of projective derangements.
\begin{theorem}
Let $d_n$ be the number of $n \times n$ projective derangements. Then
\[  1 + \sum_{n \geq 1} \frac{d_n}{\g_n}u^n = \frac{1}{1-u}\, 
               \prod_{r \geq 1} \left(1- \frac{u}{q^r} \right)^{q-1}.
\]
\end{theorem}
\proof
In this case $\Phi'=\Phi \setminus \{z-a | a \in \fq \} $. As in the proof of Theorem 8 this time we multiply and divide by the products corresponding to all the linear polynomials except $\phi(z)=z-1$. There are $q-1$ of these and so we get
\[1 + \sum_{n \geq 1} \frac{e_n}{\g_n}u^n = \prod_{\phi \neq z} \prod_{r \geq 1}
               \left(1- \frac{u^{\deg \phi}}{q^{r \deg \phi}} \right) ^{-1} \, 
               \prod_{r \geq 1} \left(1- \frac{u}{q^r} \right)^{q-1}.
\]
Use Lemma 4 to finish the proof.
\qed

It would be interesting to find a recursive formula for the $d_n$ analagous to that given in Corollary 9 for the $e_n$.

\subsection{Diagonalizable matrices}

\begin{theorem} 
Let $d_n$ be the number of diagonalizable $n \times n$ matrices. Then
\[ 1 + \sum_{n \geq 0}\frac{d_n}{\g_n}u^n = \left(  \sum_{m\geq 0} \frac{u^m}{\g_m} \right)^q  \]
and \[  d_n = \sum_{n_1 + \cdots + n_q = n} \frac{\g_n}{\g_{n_1}\cdots\g_{n_q}} . \]
\end{theorem}
\proof
We use Lemma 5. Diagonalizable matrices have conjugacy class data that only involves the linear polynomials  $\phi(z)=z-a$, $a \in \fq$ and partitions $\lambda_\phi(A)$ that are  either empty or have the form $1 \geq 1 \geq \ldots \geq 1$. These partitions are indexed by non-negative integer. For $\phi(z)=z-a$ and for $\lambda$ the partition consisting of $m$ 1's, the corresponding $m \times m$ matrix in the canonical form is the diagonal matrix with $a$ on the main diagonal. The centralizer subgroup of this matrix is the full general linear group $\gl_m(q)$ and so $\cpl = \g_m$.
Then from Lemma 5 we see that
\begin{align*}
1 + \sum_{n \geq 0}\frac{d_n}{\g_n}u^n&=\prod_{a \in \fq}\sum_{m\geq 0} \frac{u^m}{\g_m}  \\
  &= \left(  \sum_{m\geq 0} \frac{u^m}{\g_m} \right)^q .
\end{align*}  
From this the formula for $d_n$ follows immediately. 
\qed

It should be noted that the formula for $d_n$ follows directly from the knowledge of the centralizer subgroup of a diagonalizable matrix and does not require the full machinery of the cycle index generating function. Also, there is an alternative approach to this problem given in \cite{Morrison04}.

\subsection{Projections}
\begin{theorem}
Let $p_n$ be the number of $n \times n$ projections. Then
\[ 1 + \sum_{n \geq 0}\frac{p_n}{\g_n}u^n= \left(  \sum_{m\geq 0} \frac{u^m}{\g_m} \right)^2 \]
and 
\[  p_n = \sum_{0 \leq i \leq n} \frac{\g_n}{\g_i \g_{n-i}} . \]
\end{theorem}
\proof
Since projections are diagonalizable matrices with eigenvalues restricted to be in the set $\{0,1\}$, it follows that the product in the generating function for the cycle index is over the two polynomials $z$ and $z-1$ and the partitions are restricted as described in the previous theorem. Thus, we see that
\begin{align*}
1 + \sum_{n \geq 0}\frac{p_n}{\g_n}u^n&=\prod_{z,z-1}\,\sum_{m\geq 0} \frac{u^m}{\g_m}  \\
  &= \left(  \sum_{m\geq 0} \frac{u^m}{\g_m} \right)^2.
\end{align*}  
The formula for $p_n$ follows immediately.
\qed

\subsection{Solutions of $A^k=I$ }
\begin{theorem}
Let $a_n$ be the number of $n \times n$ matrices $A$ satisfying $A^2=I$, and assume that  the base field $\fq$ has characteristic two. Then
\[ a_n =  \sum_{0 \leq i \leq n/2} \frac{\g_n}{q^{i(2n-3i)}\g_{i} \g_{n-2i} } .\]
\end{theorem}
\proof
The rational canonical form for $A$ is a direct sum of companion matrices for $z-1$ and $(z-1)^2$. No other polynomials occur. Thus, the partition $\lambda_\phi (A)$ for $\phi(z)=z-1$ consists of $b_1$ repetitions of 1 and $b_2$ repetitions of 2, where $b_1 +2b_2=n$. 
From Lemma 5  we see that
\[ 1+\sum_{n \geq 1} \frac{a_n}{\g_n} u^n 
   = \sum_{b_1,\, b_2} \frac{u^{b_1+2b_2}}{c_{z-1}(b)} ,\]
where $b=(b_1,b_2)$ denotes the partition.
To compute $c_{z-1}(b)$ we refer to the formula stated in the proof of Lemma 2.  In that formula we have $d_1=b_1 + b_2$ and $d_2=b_1+2b_2$, and then
\[ c_{z-1}(b)= \prod_{i=1}^2\prod_{k=1}^{b_i}(q^{d_i}-q^{d_i -k }) . \]
This becomes
\[  c_{z-1}(b)= \prod_{k=1}^{b_1}(q^{d_1}-q^{d_1-k})\prod_{k=1}^{b_2}(q^{d_2}-q^{d_2-k}) .\]
Now let $i=d_1 - k$ in the first product and let $i=d_2-k$ in the second product, so that
\[  c_{z-1}(b)= \prod_{i=b_2}^{d_1-1} (q^{d_1}-q^i)  \prod_{i=d_2-b_2}^{d_2-1}(q^{d_2}-q^i) .
\]
From each factor of the first product we pull out a factor of $q^{b_2}$ and from each factor of the second product we remove a factor of $q^{d_2-b_2}$. This gives
\[  c_{z-1}(b)=q^{b_1 b_2}\g_{b_1} q^{(d_2-b_2)b_2} \g_{b_2} .\]
Then 
\[  \frac{a_n}{\g_n} = \sum_{b_1+2b_2=n} \frac{1}{c_{z-1}(b)} , \]
where we note that $d_2=n$ and $b_1=n-2b_2$. Therefore,
\[ \frac{a_n}{\g_n} = \sum_{b_1+2b_2=n} \frac{1}{q^{b_2(3n-2b_2)}\g_{b_1}\g_{b_2}} ,\]
which is equivalent to the formula to be proved.
\qed

\begin{theorem}
Assume that  $k$ is a positive integer relatively prime to $q$. Let 
\[ z^k-1 = \phi_1(z) \phi_2(z) \cdots \phi_r(z)  \]
be the factorization of $z^k-1$ over $\fq$ into distinct irreducible polynomials with  $d_i = \deg \phi_i $, and let $a_n$ be the number of $n \times n$ matrices over $\fq$ that are solutions of $A^k=I$. Then
\[   1 + \sum_{n \geq 1} \frac{a_n}{\g_n}u^n =
 \prod_{i=1}^r \sum_{m \geq 0} \frac{q^{md_i}}{\g_m(q^{d_i})} .\]
 \end{theorem}
 \proof
The rational canonical form of a matrix $A$ satisfying $A^k=I$ is a direct sum of any number of copies of the companion matrices of the $\phi_i$. Thus, with Lemma 5 the product is taken over the $\phi_i$ for $i=1,\ldots,r$ and the subset  $L_{\phi_i}$ is the same for all $i$ and consists of the  partitions in which all parts are 1.  Hence, $L_{\phi_i}$ can be identified with the non-negative integers $m=0,1,2,\ldots$. For the partition $\lambda$  given by $m$ 1's, the value of $c_{\phi_i}(\lambda)$ is the order of the group of automorphisms of the $\fq[z]$-module $\fq[z]/(\phi_i^{m})$, but this module is the  direct sum of $m$ copies of the extension field of degree $d_i$ over $\fq$. Therefore, the group of automorphisms is the general linear group $\gl_{m}(q^{d_i})$, and so $c_{\phi_i}(\lambda)=\g_m(q^{d_i})$.
\qed

\subsection{Cyclic matrices}
\begin{theorem}
Let $a_n$ be the number of cyclic $n  \times n$ matrices. Then
\[ 1 + \sum_{n \geq 1} \frac{a_n}{\g_n}u^n =  
     \prod_{d \geq 1}\left( 1 +  \frac{1}{q^d - 1}\,\frac{u^d}{1-(u/q)^d} \right)^{\nu_d}  \]
and
\[  1 + \sum_{n \geq 1} \frac{a_n}{\g_n}u^n = \frac{1}{1-u} \prod_{d \geq 1}
       \left(1 + \frac{u^d}{q^d(q^d-1)} \right) ^{\nu_d} .\]
\end{theorem}
\proof
In terms of the conjugacy class data a matrix is cyclic if $\lambda_\phi$ has at most one part for each $\phi$. Thus, $L_\phi= \{ \emptyset, 1, 2, \ldots \}$, where $m$ means the partition of $m$ having just one part. For these partitions we have $c_\phi(\emptyset) =1$ and
 $c_\phi(m)= q^{m \deg \phi} - q^{(m-1)\deg \phi} $. Let $a_n$ be the number of cyclic matrices of size $n \times n$. Beginning with Lemma 5 we find the generating function.
 \begin{align*}
1 + \sum_{n \geq 1} \frac{a_n}{\g_n}u^n &= \prod_{\phi} \sum_{L_\phi} \frac{u^{|\lambda| \deg \phi}}{\cpl}  \\
  &=  \prod_{\phi}\left( 1 +  \sum_{m \geq 1} \frac{u^{m\deg \phi}}{q^{m \deg \phi} - q^{(m-1)\deg \phi} } 
     \right) \\
 &= \prod_{d \geq 1}\left( 1 +  \sum_{m \geq 1} \frac{u^{md}}{q^{m d} - q^{(m-1)d} } 
     \right) ^{\nu_d} \\
  &= \prod_{d \geq 1}\left( 1 +  \frac{1}{q^d - 1} \sum_{m \geq 1} \frac{u^{md}}{q^{(m-1)d} } 
     \right) ^{\nu_d} \\
   &=  \prod_{d \geq 1}\left( 1 +  \frac{1}{q^d - 1}\,\frac{u^d}{1-(u/q)^d} \right)^{\nu_d}.
\end{align*}
For the second part, we use Lemma 6 to get
\[ 1 + \sum_{n \geq 1} \frac{a_n}{\g_n}u^n = \frac{1}{1-u}
          \prod_{d \geq 1} \left(1-\frac{u^d}{q^d} \right)^{\nu_d}
           \prod_{d \geq 1}\left( 1 +  \frac{1}{q^d - 1}\,\frac{u^d}{1-(u/q)^d} \right)^{\nu_d}.
\]
Combine the products and simplify to finish the proof.
\qed
\subsection{Semi-simple matrices}
\begin{theorem}
Let $a_n$ be the number of $n \times n$ semi-simple matrices over $\fq$. Then
\[ 
   1 +  \sum_{n \geq 1} \frac{a_n}{\g_n} u^n=  
   \prod_{d \geq 1}\left( 1 + \sum_{j \geq 1} \frac{u^{jd}}{\g_j(q^d) }\right)^{\nu_d}.
\]
\end{theorem}
\proof  In order for a matrix to be semi-simple the partition associated to a polynomial $\phi$ has no parts greater than 1. Thus, for all $\phi$, $L_\phi = \{ \emptyset, 1, 1^2, \ldots, 1^j, \ldots \}$, where $1^j$ means the partition of $j$ consisting of $j$ copies of 1. Then $c_\phi ( 1^j)$ is the order of the automorphism group of the $\fq[z]$-module which is the direct sum of $j$ copies of $\fq[z]/(\phi)$, which can be identified with $(\F_{q^d} )^ j$, where $d=\deg \phi$. Furthermore, the automorphisms of the module $\left( \fq[z]/(\phi)\right) ^{\oplus j}$ can be identified with the $\F_{q^d}$-linear automorphisms of  $(\F_{q^d} )^ j$. This group of automorphisms is $\gl_j (q^d)$, and so $c_\phi ( 1^j)=\g_j(q^d)$.

From Lemma 5 we find
 \begin{align*}
1 + \sum_{n \geq 1} \frac{a_n}{\g_n}u^n &= \prod_{\phi} \sum_{L_\phi} \frac{u^{|\lambda| \deg \phi}}{\cpl}  \\
     &= \prod_\phi \left(1+\sum_{j \geq 1}\frac{u^{j \deg \phi}}{\g_j(q^d)} \right) \\
     &= \prod_{d \geq 1} \left(1+\sum_{j \geq 1}\frac{u^{jd}}{\g_j(q^d)} \right)^{\nu_d} .
\end{align*}
\qed
\subsection{Separable matrices}
\begin{theorem}
 Let $a_n$ be the number of separable $n \times n$ matrices over $\fq$. Then 
 \[ 
    1 + \sum_{n \geq 1} \frac{a_n}{\g_n} u^n =
    \prod_{d \geq 1} \left( 1 + \frac{u^d}{q^d-1} \right)^{\nu_d}
\]     
and
\[  
 1 + \sum_{n \geq 1} \frac{a_n}{\g_n} u^n =  \frac{1}{1-u}
    \prod_{d \geq 1} \left( 1 + \frac{u^d(u^d -1)}{q^d(q^d-1)} \right)^{\nu_d}.
\]     
\end{theorem}
\proof  For a matrix to be separable the allowed partitions are either empty or the unique partition of 1. So the sum over $L_\phi$ in Lemma 5 is a sum of two terms. For $\lambda = \emptyset$, 
\[ \frac{u^{| \lambda| \deg \phi}}{\cpl}=1, \] and for $\lambda = 1$, 
\[ \frac{u^{| \lambda| \deg \phi}}{\cpl}= \frac{u^{\deg \phi}}{q^{\deg \phi} - 1} .\]
Therefore,
\begin{align*}
 1 + \sum_{n \geq 1} \frac{a_n}{\g_n}u^n & =  \prod_{\phi} \sum_{L_\phi} \frac{u^{|\lambda| \deg \phi}}{\cpl}  \\
    &= \prod_{\phi} \left( 1 + \frac{u^{\deg \phi}}{q^{\deg \phi}-1} \right)  \\
    &= \prod_{d \geq 1} \left( 1 +\frac{u^{d}}{q^{d}-1} \right)^{\nu_d}.
\end{align*}
For the second statement of the theorem multiply the right side in the line above by 
 \[  \frac{1}{1-u}\prod_{d \geq 1} \left( 1 - \frac{u^d}{q^d} \right)^{\nu_d} ,\]
 which is equal to 1 by Lemma 6. Then combine the products.
\qed 

\subsection{Conjugacy classes}
\begin{theorem}
Let $a_n$ be the number of conjugacy classes of $n \times n$ matrices. Then the ordinary generating function for the sequence $\{ a_n \}$ is given by
\[  1 + \sum_{n \geq 1} a_n u^n = \prod_{r \geq 1} \frac{1}{1-qu^r}  .\]
\end{theorem}
\proof 
Recall the ordinary power series generating function for partitions factors as the infinite product
\[  \sum_{n \geq 0} p_n u^n = \prod_{r \geq 1} \frac{1}{1-u^r} , \]
where $p_n$ is the number of partitions of the integer $n$. A conjugacy class of $n \times n$ matrices is uniquely specified by the choice of a partition $\lambda_\phi$ for each monic irreducible polynomial $\phi$ such that $\sum_{\phi} | \lambda_\phi | \deg \phi = n$.  Consider the infinite product over $\phi$
\[  \prod_{\phi} \sum_{n \geq 0} p_n u^{n \deg \phi}  .\]  The coefficient of $u^n$ is a sum of terms of the form $ p_{n_1} u^{n_1 \deg \phi_1} \cdots p_{n_k} u^{n_k \deg \phi_k} $ where $n = \sum n_k \deg \phi_k$. Therefore, the $u^n$-coefficient is the number of conjugacy classes of $n \times n$ matrices.
Then
\begin{align*}
  1 + \sum_{n \geq 1} a_n u^n  &= \prod_{\phi} \sum_{n \geq 0} p_n u^{n \deg \phi} \\
      &=  \prod_{\phi} \prod_{r \geq 1} \frac{1}{1-u^{r \deg \phi}}  \\
      &= \prod_{r \geq 1} \prod_{\phi} \frac{1}{1-u^{r \deg \phi} }  \\
      &= \prod_{r \geq 1} \prod_{d \geq 1} \left( \frac{1}{1-u^{r d}} \right)^{\nu_d}
\end{align*}
Then by substituting $qu^r$ for $u$ in Lemma 6 and inverting we see that
\[  \prod_{d \geq 1} \left( \frac{1}{1-u^{r d}} \right)^{\nu_d}= \frac{1}{1-qu^r} .\]
Therefore
\[   1 + \sum_{n \geq 1} a_n u^n=  \prod_{r \geq 1}  \frac{1}{1-qu^r} .\]
\qed

\raggedright
\emph{Sequences considered:} A002416, A002820,  A002884,  A005329, A006116 --A006122, A015195--A015217, A022166--A022188, A053718, A053722, A053725, A053770--A053777, A053763, A053846--A053849, A053851--A053857, A053859--A053863,  A069777, A070933, A076831, A082877.

 \end{document}